\documentclass{amsart}
\usepackage{amsfonts}
\usepackage{amsmath,amssymb}
\usepackage{amsthm}
\usepackage{amscd}
\usepackage{graphics}
\usepackage{graphicx}

\theoremstyle{remark}{
\newtheorem{Def}{{\rm Definition}}

}
\theoremstyle{plain}
{

\newtheorem{Thm}{Theorem}

}

\begin{document}
\title[.]{Realizations of planar graphs as Poincar\'e-Reeb graphs of refined algebraic domains}
\author{Naoki kitazawa}
\keywords{(Non-singular) real algebraic manifolds and real algebraic maps. Algebraic domains. Poincar\'e-Reeb Graphs. Singularity theory of Morse(-Bott) functions.%E3%81%8A%E3%82%80%E3%81%99%E3%81%B3-%E5%89%8D%E4%BB%A3%E6%9C%AA%E8%81%9E%E3%81%AE-%E7%95%B0%E5%B8%B8%E4%BA%8B%E6%85%8B-%E7%99%BA%E7%94%9F-%E3%83%92%E3%83%AD%E3%82%A4%E3%83%B3%E3%81%AE%E5%90%8D%E5%89%8D%E3%81%8C%E6%B6%88%E3%81%88%E3%81%9F/ar-AA1xI5nh?ocid=msedgntp&pc=ASTS&cvid=0a822926a6a4455189d491b82ddf8032&ei=6Reeb graphs. Circles in plane geometry. \\
\indent {\it \textup{2020} Mathematics Subject Classification}: Primary~14P05, 14P10, 52C15, 57R45. Secondary~ 58C05.}

\address{Institute of Mathematics for Industry, Kyushu University, 744 Motooka, Nishi-ku Fukuoka 819-0395, Japan\\
 TEL (Office): +81-92-802-4402 \\
 FAX (Office): +81-92-802-4405 \\
}
\email{n-kitazawa@imi.kyushu-u.ac.jp, naokikitazawa.formath@gmail.com}
\urladdr{https://naokikitazawa.github.io/NaokiKitazawa.html}
\maketitle
\begin{abstract}
{\it Algebraic domains} are regions in the plane surrounded by mutually disjoint non-singular real algebraic curves.
{\it Poincar\'e-Reeb Graphs} of them are graphs they naturally collapse: such graphs are formally formulated by Sorea, for example, around 2020. Their studies found that nicely embedded planar graphs are Poincar\'e-Reeb graphs of some algebraic domains. These graphs are generic with respect to the projection to the horizontal axis. Problems, methods and results are elementary and natural and they apply natural approximations nicely for example. 

We present our new approach to extension of the result to a non-generic case and an answer. 
We first formulate generalized algebraic domains, surrounded by non-singular real algebraic curves which may intersect with normal crossings. Such domains and certain classes of them appear in related studies of graphs and regions surrounded by algebraic curves explicitly.

\end{abstract}
%【REVISE】 combinatoric ～ is → combinatorial object. It is .
%【REVISE】  such that a point is a vertex if and only if the corresponding connected component of the level set contains some singular points → whose vertex set is the set of all points containing some singular points in the corresponding connected component of the level set .
%【REVISE】 We delete "extending the result before".
\section{Introduction.}
\label{sec:1}
In real algebraic geometry, regions in the plane surrounded by (so-called) {\it non-singular} real algebraic curves are fundamental spaces and objects. \cite{bodinpopescupampusorea, sorea1, sorea2} show a kind of studies which are also elementary, natural and surprisingly, developing recently. They try to understand the shapes, especially, convexity, of the regions. They are defined as {\it algebraic domains}. Graphs they naturally collapse to respecting the projection to the horizontal axis $\{(t,0) \mid t \in \mathbb{R}\}$ are introduced and shown to be important: hereafter let ${\mathbb{R}}^n$ denote the $n$-dimensional Euclidean space, which is also a smooth manifold equipped with the standard Euclidean metric, and $||x|| \geq 0$ the distance between $x \in {\mathbb{R}}^n$ and the origin $0 \in {\mathbb{R}}^n$ under the metric as usual ($\mathbb{R}:={\mathbb{R}}^1$). 
They have shown that for a naturally embedded planar graph being generic with respect to the projection, we can find an algebraic domain collapsing naturally to the graph. 
There classical and strong arguments such as approximation by real polynomials are essential. Such a graph is also named a {\it Poincar\'e-Reeb graph} of the algebraic domain. 
Our study extends their result to graphs which may not be generic in the sense above.

\subsection{Our notation on topological spaces, manifolds and graphs.}
Let ${\pi}_{m,n}:{\mathbb{R}}^m \rightarrow {\mathbb{R}}^n$ with $m>n \geq 1$ denote the so-called canonical projection ${\pi}_{m,n}(x)=x_1$ where $x=(x_1,x_2) \in {\mathbb{R}}^n \times {\mathbb{R}}^{m-n}={\mathbb{R}}^m$. We also use $D^k:=\{x \in {\mathbb{R}}^k \mid ||x|| \leq 1\}$, for the $k$-dimensional unit disk, and $S^k:=\{x \in {\mathbb{R}}^{k+1} \mid ||x||=1\}$, for the $k$-dimensional unit sphere, for example. 

For a topological space $X$ and its subspace $Y \subset X$, we use $\overline{Y}$ for its closure and $Y^{\circ}$ for its interior: we omit information on the outer space $X$ unless otherwise stated: we can guess from our arguments. For a topological space $X$ decomposed into a so-called cell complex, we can define the dimension $\dim X$ uniquely as the dimension of the cell of the maximal dimension (only depending on the topology of $X$): topological manifolds, polyhedra, and graphs, which are regarded as $1$-dimensional CW complexes, are of such a class. For a topological manifold $X$ whose boundary is non-empty, we use $\partial X$ for its boundary and ${\rm Int}\ X:=X-\partial X$. For a smooth manifold $X$ and $x \in X$, we use $T_xX$ for the tangent vector space of $X$ at $x$. For smooth manifolds $X$ and $Y$ and a smooth map $c:X \rightarrow Y$, a point $x \in X$ is its singular point if the rank of the differential ${dc}_x:T_xX \rightarrow T_{c(x)}Y$ is smaller than both the dimensions $\dim X$ and $\dim Y$: note that ${dc}_x$ is a linear map.
The zero set $S$ of a real polynomial map or more generally, a union $S$ of its connected components is {\it non-singular} if the polynomial map has no singular point in the set $S$: remember the implicit function theorem.  

A graph is a CW complex where an {\it edge} is a $1$-cell and a {\it vertex} is a $0$-cell. The set of all edges (vertices) of the graph is the {\it edge set} ({\it vertex set}) of the graph. Two graphs $G_1$ and $G_2$ are {\it isomorphic} if there exists a (piecewise smooth) homeomorphism $\phi:G_1 \rightarrow G_2$ mapping the vertex set of $G_1$ onto that of $G_2$: such a homeomorphism is called an {\it isomorphism} of the graphs. A {\it digraph} is a graph all of whose edges are oriented and two digraphs are {\it isomorphic} if there exists an isomorphism of graphs between them preserving the orientations, which is defined as an {\it isomorphism} of the digraphs.
\subsection{Refined algebraic domains.}
\label{subsec:1.2}
%We introduce {\it arrangements of circles for Morse-Bott functions} or {\it MBC arrangements}.
\begin{Def}
\label{def:1}
A pair of a family $\mathcal{S}=\{S_j \subset {\mathbb{R}}^2\}$ each $S_j$ of which is the zero set of a real polynomial $f_j$ and non-singular and to each of which $f_j$ is associated and a region $D_{\mathcal{S}} \subset {\mathbb{R}}^2$ satisfying the following conditions is called a {\it refined algebraic domain}. 
\begin{enumerate}
\item The region satisfies $D_{\mathcal{S}}={\bigcap}_{S_j \in \mathcal{S}} \{x \mid f_j(x)>0\} \subset {\mathbb{R}}^2$ and a bounded connected component of ${\mathbb{R}}^2-{\bigcup}_{S_j \in \mathcal{S}} S_j$ and the intersection $\overline{D_{\mathcal{S}}} \bigcap S_j$ is non-empty for any curve $S_j \in \mathcal{S}$.
\item At points in $\overline{D_{\mathcal{S}}}$, at most two distinct curves $S_{j_1}, S_{j_2} \in \mathcal{S}$ intersect and the following are satisfied: for each point $p_{j_1,j_2}$ in such an intersection, the sum of the tangent vector spaces of them at $p_{j_1,j_2}$ coincides with the tangent vector space of ${\mathbb{R}}^2$ at $p_{j_1,j_2}$. %and three distinct curves do not intersect in $\overline{D_{\mathcal{S}}}$.
\end{enumerate}
\end{Def}

This also respects \cite{kohnpieneranestadrydellshapirosinnsoreatelen} for example. We discuss the restriction of ${\pi}_{2,1}$ to $\overline{D_{\mathcal{S}}}$. We consider the set $F_{D_{\mathcal{S}}}$ of all points in the following. This is finite thanks to the real algebraic situation.
\begin{itemize}
\item 
Points in $\overline{D_{\mathcal{S}}}$ which are also in exactly two distinct curves $S_{j_1}$ and $S_{j_2}$.
\item If we remove the finite set before from the set $\overline{D_{\mathcal{S}}}-D_{\mathcal{S}}$ of dimension $1$, then we have a smooth manifold of dimension $1$ (a curve which is not necessarily connected) and which has no boundary. Points which are singular points of the restriction of ${\pi}_{2,1}$ to the obtained smooth curve in $\overline{D_{\mathcal{S}}}-D_{\mathcal{S}}$. 
\end{itemize} 
We can define the following equivalence relation ${\sim}_{D_{\mathcal{S}}}$ on $\overline{D_{\mathcal{S}}}$: two points are equivalent if and only if they belong to a same component of the preimage of a single point for the restriction of ${\pi}_{2,1}$ to $\overline{D_{\mathcal{S}}}$. Let $q_{D_{\mathcal{S}}}$ denote the quotient map and $V_{D_{\mathcal{S}}}$ the function uniquely defined by the relation ${\pi}_{2,1}=V_{D_{\mathcal{S}}} \circ q_{D_{\mathcal{S}}}$. The quotient space $W_{D_{\mathcal{S}}}:=\overline{D_{\mathcal{S}}}/{\sim}_{D_{\mathcal{S}}}$ is a digraph by the following. We can check this from general theory \cite{saeki1, saeki2} or see \cite{kitazawa5} for example: we do not need to understand this theory.
\begin{enumerate}
\item The vertex set is the set of all points $v$ whose preimage ${q_{D_{\mathcal{S}}}}^{-1}(v)$ contains at least one point of the finite set $F_{D_{\mathcal{S}}}$ above.
\item The edge connecting $v_1$ and $v_2$ are oriented as one departing from $v_1$ and entering $v_2$ according to $V_{D_{\mathcal{S}}}(v_1)<V_{D_{\mathcal{S}}}(v_2)$.
\end{enumerate}
\begin{Def}
We call the (di)graph above a {\it Poincar\'e-Reeb }({\it di}){\it graph of $D_{\mathcal{S}}$}.
\end{Def}
As this graph, we can consider a situation where for a graph $G$, a nice map $V_{G}$ on its vertex set onto a partially ordered set $P$ is given and orients the graph according to the values. More precisely, each edge $e$ of the graph connects two distinct vertices $v_{e,1}$ and $v_{e,2}$ and it is oriented. Furthermore, it is oriented according to the rule: the edge $e$ departs from $v_{e,1}$ and enters $v_{e,2}$ if $V_{G}(v_{e,1})<V_{G}(v_{e,2})$: let $<$ denote the order on $P$. 
We call a pair of such a graph $G$ and a map $V_{G}$ a V-digraph. For V-graphs, {\it isomorphisms} between two V-digraphs and the notion that two V-digraphs are {\it isomorphic} can be defined, based on the property of preserving the orders of the values of the maps $V_{G}$.
We can also define the {\it Poincar\'e-Reeb V-digraph of $D_{\mathcal{S}}$} by associating the function $V_{D_{\mathcal{S}}}$.
\subsection{Our main result.}
Two graphs, digraphs, and V-digraphs are {\it weakly isomorphic} if there exists a homeomorphism regarded as an isomorphism after suitable addition of finitely many vertices: the edge sets of the graphs also change. 
\begin{Thm}
\label{thm:1}
For any finite and connected graph $G$ and a piecewise smooth function $c_G:G \rightarrow \mathbb{R}$ such that the restriction ${c_G} {\mid}_{e}$ is injective for each edge $e$ of $G$, we can canonically give $G$ the structure of a V-digraph by the function $c_G$. We also assume the following.
\begin{enumerate}
\item The function $c_G$ is the composition of some piecewise smooth embedding $e_G:G \rightarrow {\mathbb{R}}^2$ with ${\pi}_{2,1}$.
\item The degree of each vertex of $G$ is not $2$. The local extremum of $c_G$ must be achieved at a vertex of degree $1$. 
\end{enumerate}
Then we have a refined algebraic domain $D_{G}$ and its Poincar\'e-Reeb V-digraph of $D_{G}$ and the V-digraph $G$ are weakly isomorphic. 
\end{Thm}
Note that \cite{bodinpopescupampusorea} has shown a generic case: the degrees of vertices are always $1$ or $3$ with the values of $c_G$ at distinct vertices being always distinct. 
They only consider real algebraic domains: curves are mutually disjoint. On the other hand, the resulting V-digraphs have been shown to be isomorphic in the case. The constraint that curves are the zero sets of some real polynomials is not considered there and the author has commented first in \cite{kitazawa4}: in the original study the curves are only unions of some connected components of the zero sets.

\subsection{Organization of our paper and our main work.}

In the next section, we show Theorem \ref{thm:1}. 
In the third section, we introduce that our graph is regarded as the so-called {\it Reeb graph} of a nice real algebraic function. The Reeb graph of a smooth function is a classical and fundamental object (\cite{reeb}). This is defined as the quotient space of the manifold similarly to Poincar\'e-Reeb graphs. This represents the manifold compactly.
This also gives a new answer to the following: can we reconstruct a real algebraic function whose Reeb graph is isomorphic to the given graph? We also present related studies since the birth of the study by Sharko (\cite{sharko}), in 2006, reconstructing nice differentiable (smooth) functions on closed surfaces.

%Hereafter, the Euclidean space ${\mathbb{R}}^k$ is regarded as the vector space canonically.

\section{A proof of Theorem \ref{thm:1}.}
\label{sec:2}
In this section, we prove Theorem \ref{thm:1}, our main new result.

We use fundamental arguments from singularity theory and real algebraic geometry. Especially, approximations.

See \cite{golubitskyguillemin} for singularity theory of differentiable maps. For example, we mainly consider the Hesse matrix of a differentiable maps of the class $C^2$, the symmetric matrix canonically defined as the matrix of the second derivatives. Smooth functions such that the determinant of the Hesse matrix, or the {\it Hessian}, at each point of the space of the domain is not $0$, are important. Such a function is also a so-called {\it Morse} function. A {\it Morse-Bott} function is a smooth function whose singular point is represented as the composition of a smooth map with no singular point with a Morse function for suitable local coordinates.

See \cite{bochnakcosteroy, kollar, lellis} for real algebraic geometry, for example.  See also \cite{elredge} for approximations by real polynomials. 

Last, our present approximation mainly respects \cite{bodinpopescupampusorea} as an explicit and important case and revises some. 
We respect these arguments from singularity theory and approximations implicitly.

Related to this, see also \cite{kitazawa4} where we do not assume the arguments from the preprint.

Hereafter, an {\it ellipsoid of ${\mathbb{R}}^2$ centered at a point $x_0=(x_{0,1},x_{0,2})$} means a set of the form $\{x=(x_1,x_2) \in {\mathbb{R}}^2 \mid a_1{(x_1-x_{0,1})}^2+a_2{(x_2-x_{0,2})}^2 \leq r\}$ where $a_1,a_2,r>0$. Sets of this type are also important.

\begin{proof}[A proof of Theorem \ref{thm:1}]
We consider the graph $e_G(G) \subset {\mathbb{R}}^2$.
We can change the graph $e_G(G)$ which is also a CW complex as follows.

Here, we choose sufficiently small positive numbers ${\epsilon}_1, {\epsilon}_2>0$. We can choose them so that we can argue with no problem. We can see this by following our arguments.

First we consider a point $p \in \mathbb{R}$ where the preimage ${{\pi}_{2,1}}^{-1}(p)$ contains at least one vertex of $e_G(G)$.
We can choose vertices $v_{p,{\rm m}}=(p,p_{{\rm m}})$ and $(p,p_{{\rm M}})$ contained in the preimage in such a way that the value $p_
{{\rm m}}$ is the minimum of the values of the second components among such vertices in the preimage and that the value $p_{{\rm M}}$ is the maximum of the values of the second components among such vertices in the preimage.
We first add a segment $S_p:=\{(p,y) \mid p_{{\rm m}}-{\epsilon}_2 \leq y \leq p_{{\rm M}}+{\epsilon}_2\}$.

We have a new CW complex $e_G(G) \bigcup S_p$. 
Let $v_p:=(p,p_v) \in e_G(G) \subset {\mathbb{R}}^2$ be a vertex such that $c_G$ does not have a local extremum at ${e_G}^{-1}(v_p) \in G$, which is also regarded as a vertex of $G$.

Let $N_v$ be a neighborhood represented as $N_v:=\{x=(x_1,x_2) \mid p-{\epsilon}_1 < x_1 < p+{\epsilon}_1, p_v-{\epsilon}_2 < x_2 < p_v+{\epsilon}_2\}$.

We define $S_{p,{\epsilon}^{\prime}}:=\{(p,y) \mid p_{{\rm m}}+ {\epsilon}^{\prime} \leq y \leq p_{{\rm M}}-{\epsilon}^{\prime}\}$ and $S_{p \pm \frac{{\epsilon}_1}{2},{\epsilon}^{\prime}}:=\{(p \pm \frac{{\epsilon}_1}{2},y) \mid p_{{\rm m}}+ {\epsilon}^{\prime} \leq y \leq p_{{\rm M}}-{\epsilon}^{\prime}\}$ where another sufficiently small positive number ${\epsilon}^{\prime}>0$ is chosen.

The set ${{\pi}_{2,1}}^{-1}(p \pm {\epsilon}_1) \bigcap e_G(G)$ is finite. There exists a curve in $e_G(G)$ and an edge departing from each point $p_{\pm {\epsilon}_1,j}$ there to the vertex $v_p$. We change each of these curves to a union of two segments intersecting in a one-point set. 
More precisely, we can change and change the curves as follows. 

For each curve in $e_G(G)$, the first segment is a straight segment departing from the point $p_{\pm {\epsilon}_1,j}$ and entering another point $p_{\pm \frac{{\epsilon}_1}{2},j}$, which is a point in $N_v \bigcap S_{p \pm \frac{{\epsilon}_1}{2},{\epsilon}^{\prime}}$. We also choose these segments as mutually disjoint ones. 

Next, for each curve, the second segment is chosen as the unique horizontal segment connecting $p_{\pm \frac{{\epsilon}_1}{2},j}$ and a point in the segment $S_{p,{\epsilon}^{\prime}}$, which can be uniquely chosen. We also construct the segments in such a way that the horizontal segments from points $p_{-\frac{{\epsilon}_1}{2},j}$, the values of whose first components are $p-\frac{{\epsilon}_1}{2}$, are beyond the horizontal segments from points $p_{+\frac{{\epsilon}_1}{2},j}$, the values of whose first components are $p+\frac{{\epsilon}_1}{2}$.
Thus $N_v \bigcap (e_G(G) \bigcup S_p)$ is changed.

For all of such vertices $v_p$ and values $p \in \mathbb{R}$ with the preimages ${{\pi}_{2,1}}^{-1}(p)$ containing some vertices in the graph $e_G(G)$, we can do similarly and do. 
We also remove $S_p-{\bigcup}_{v_p} \overline{N_v}$ for all the values $p$ here. Thus we have a new 1-dimensional CW complex $G_{\epsilon}$ from $e_G(G)$.

For $G_{\epsilon}$, we can consider a sufficiently small regular neighborhood \cite{hirsch} as a $2$-dimensional smooth compact submanifold $M_D$ in ${\mathbb{R}}^2$.
We can also choose $M_D$ in such a way that the boundary $M_S$ is the zero set of some real polynomial and a non-singular set according to the presented theory on approximation. More precisely, we can do so that the restriction of ${\pi}_{2,1}$ to the boundary $M_S$ satisfies the following.
\begin{itemize}
\item The restriction is a Morse function.
\item All of singular points of the restriction is as follows.
\begin{itemize}
\item Each singular point of the function is corresponded to either a vertex of $e_G(G)$ of degree $1$ which is also sufficiently close or a connected component of $N_v-G_{\epsilon}$ whose closure contains no point
 of the form $(x,p_v \pm {\epsilon}_2)$ where all vertices $v_p$ as argued above are considered.
\item In addition, the correspondence is a one-to-one correspondence. We name the singular point of the function of the first (second) type a {\it definite} (resp. {\it indefinite}) type.
\item At a definite type singular point $s_v \in M_S$ of the function the value is greater (resp. smaller) than the value $c_G({\pi}_{2,1}(v))$ at the corresponding vertex $v$ of $e_G(G)$ if $c_G$ has a local maximum (resp. minimum). Note that the values are also sufficiently close.
\item Each indefinite type singular point $s_{v,j} \in M_S$ of the function is in the corresponding connected component of $N_v-G_{\epsilon}$ whose closure contains no point
 of the form $(x,p_v \pm {\epsilon}_2)$ and $s_{v,j}$ is sufficiently close to the segments $S_p$ and $S_{p,{\epsilon}^{\prime}}$.
\end{itemize}
\end{itemize}
We can put a sufficiently small suitable ellipsoid centered at an indefinite type singular point $s_{{\rm i}}:=(s_{{\rm i},1},s_{{\rm i},2})$ of the function in ${\mathbb{R}}^2$ in such a way that the boundary of the ellipsoid contains a point $(p,s_{{\rm i},2})$. 
At a definite type singular point $s_{{\rm d}}:=(s_{{\rm d},1},s_{{\rm d},2})$ of the function, we can put a sufficiently small circle containing exactly two points in $M_S$ in such a way that one of the points is of the form $(p,s_{{\rm d},p})$ and sufficiently close to $s_{{\rm d}}$ and that the restriction of the projection ${\pi}_{2,1}$ to the intersection of the small circle and $M_D$ is injective.

We have a new set by removing the intersection of $M_D$ and each ellipsoid and the disk bounded by each new small circle from $M_D$. We naturally have a new refined algebraic domain $D_G$ and we can check that this is our desired refined algebraic domain. 

For arguments here, see also FIGUREs \ref{fig:1}--\ref{fig:3}.
\begin{figure}
	\includegraphics[width=80mm,height=80mm]{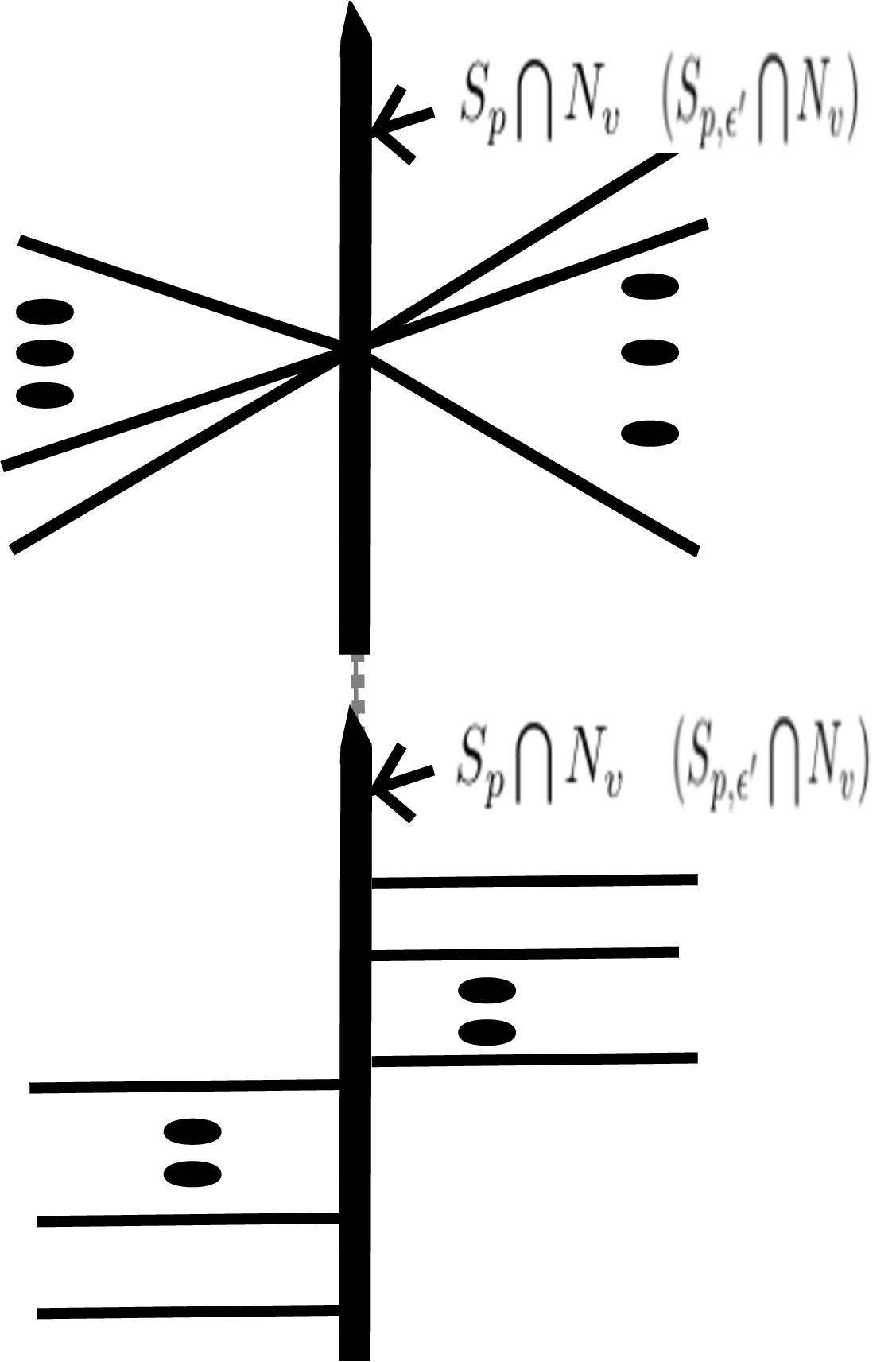}
	\caption{Around a vertex $v_p \in e_G(G)$ and $N_v$. The set $N_v \bigcap (e_G(G) \bigcup S_p)$ (, in the upper part,) is changed.}
	\label{fig:1}
\end{figure}
\begin{figure}
	\includegraphics[width=80mm,height=80mm]{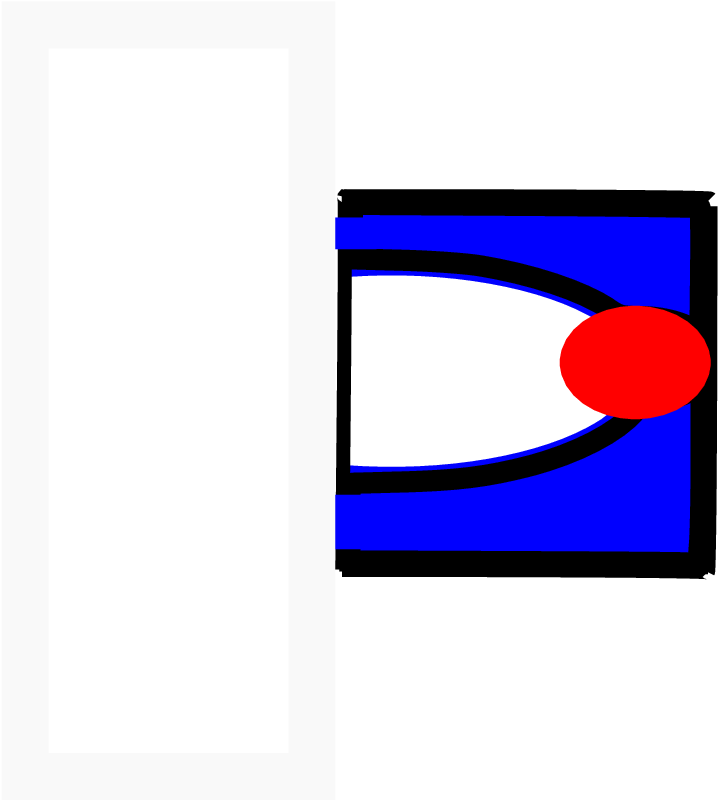}
	\caption{Around a connected component of $N_v-G_{\epsilon}$ whose closure contains no point of the form $(x,p_v \pm {\epsilon}_2)$. The blue region shows (the interior of) $M_D$ partially and the red ellipsoid is added.}
	\label{fig:2}
\end{figure}
\begin{figure}
	\includegraphics[width=80mm,height=80mm]{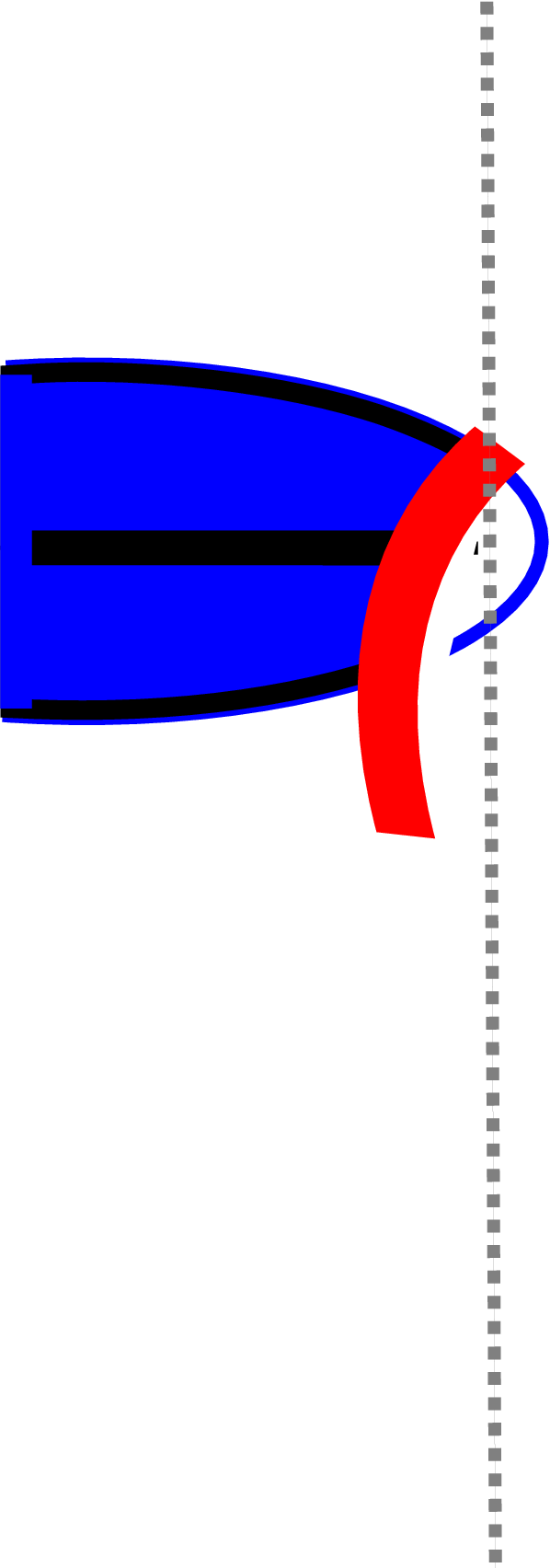}
	\caption{Around a vertex of degree $1$ of $e_G(G)$. The blue region shows (the interior of) $M_D$, containing the graph $e_G(G)$, partially. The red circle (with the disk bounded by this) is added.}
	\label{fig:3}
\end{figure}

This completes the proof.
\end{proof}
\section{Relations with Reeb graphs of real algebraic functions.}
\label{sec:3}
This section presents a kind of applications to construction of examples in real algebraic geometry and singularity theory of differentiable (smooth) maps. We assume several arguments from the published article \cite{kitazawa2} and our preprints \cite{kitazawa3, kitazawa4, kitazawa5}.

We can define the {\it Reeb} ({\it V-di}){\it graph of a smooth function $c:X \rightarrow \mathbb{R}$} on a closed manifold as follows.
These graphs have been classical and strong tools in understanding the manifolds (\cite{reeb}).

Two points of $X$ are equivalent if and only if they are in a same connected component of the preimage $c^{-1}(y)$. Let ${\sim}_c$ denote the equivalence relation and $W_c:=X/{\sim}_c$ the quotient space. Let $q_c:X \rightarrow W_c$ denote the quotient map associated with the unique continuous function $\bar{c}:W_c \rightarrow \mathbb{R}$ satisfying $c=\bar{c} \circ q_c$.
The vertex set of $W_c$ can be defined as the set of all elements $v$ such that the preimage ${q_c}^{-1}(v)$ contains at least one singular point of $c$ in the case where the image of the set of all singular points of $c$ is a finite set (\cite{saeki1, saeki2}). The graph is the {\it Reeb graph of $c$} and the pair of the graph with (the restriction of) the function $\bar{c}$ (to the vertex set) is the {\it Reeb V-digraph of $c$}.

We define a surjective map $m_{\mathcal{S},A}:\mathcal{S} \rightarrow A$ onto some finite set $A$. 
We also pose the constraint that for two distinct curves $S_{j_1}, S_{j_2} \in \mathcal{S}$ which intersect in $\overline{D_{\mathcal{S}}}$, the values of the map are distinct. We also define another positive integer valued function $m_{\mathcal{S},A,0}$ on $A$. Let
$S:=\{(x,(y_a)_{a \in A}) \in {\mathbb{R}}^2 \times  {\prod}_{a \in A} {\mathbb{R}}^{m_{\mathcal{S},A,0}(a)+1} \mid  {\prod}_{j \in {m_{\mathcal{S},A}}^{-1}(a)} (f_j(x))-{\Sigma}_{j=1}^{m_{\mathcal{S},A,0}(a)+1} {y_{a,j}}^2=0, a \in A\}$: here for the notation $y_{a,j}$ is the $j$-th component of $y_a \in {\mathbb{R}}^{m_{\mathcal{S},A,0}(a)+1}$ for example.
This is the zero set of a real polynomial map in ${\mathbb{R}}^2 \times  {\prod}_{a \in A} {\mathbb{R}}^{m_{\mathcal{S},A,0}(a)+1}$ and non-singular. We omit precise arguments. See \cite{kitazawa3}. \cite{kitazawa5} also presents this.
We consider the restriction of ${\pi}_{{\Sigma}_{a \in A}(m_{\mathcal{S},A,0}(a)+1)+2,1}$ to $S$. The Reeb V-digraph of the resulting function is isomorphic to the Poincar\'e-Reeb V-digraph of $D_{\mathcal{S}}$. See \cite{kitazawa3, kitazawa5} again and see also \cite{kitazawa2} and the preprint \cite{kitazawa4}. By our construction, we can check that the resulting function is a Morse-Bott function.

Related to this we explain history of reconstruction of nice smooth functions and the manifolds from given graphs.

\cite{sharko} is a pioneering study, constructing nice smooth functions whose Reeb graphs are isomorphic to given finite graphs of a certain nice class. The functions are locally elementary polynomials. This is extended in \cite{masumotosaeki} to the case of all finite graphs. \cite{michalak1} studies the Morse function case mainly the case of functions on closed surfaces. \cite{marzantowiczmichalak, michalak2} study a kind of general theory of Morse functions and their Reeb graphs. \cite{michalak2} mainly studies deformations of Morse functions via deformations of Reeb graphs. Following \cite{michalak2}, \cite{marzantowiczmichalak} studies classifications of Morse functions on manifolds of general dimensions via systems of hypersurfaces, represented as preimages of the functions, for example. \cite{gelbukh1, gelbukh2} study the Morse-Bott function case, mainly the case of functions on closed surfaces. \cite{kitazawa1} studies cases of functions of certain classes naturally generalizing the classes of Morse-Bott functions on $3$-dimensional closed manifolds where surfaces of preimages of points are prescribed. This is regarded as a pioneering study considering cases where preimages of single points are prescribed before reconstruction of functions. 

The case of globally real algebraic functions is pioneered in \cite{kitazawa2}. 

Our theorem gives a kind of new answers to the real algebraic case. We can reconstruct real algebraic functions whose Reeb graphs are only homeomorphic to the given graphs. Related to this, \cite{gelbukh1} is for reconstruction of Morse-Bott functions in the differentiable (smooth) situation whose Reeb graphs are only homeomorphic to the given graphs.
\section{Conflict of interest and Data availability.}
\noindent {\bf Conflict of interest.}
The author works at Institute of Mathematics for Industry (https://www.jgmi.kyushu-u.ac.jp/en/about/young-mentors/). This is closely related to our study. We thank them for supports and encouragement. The author is also a researcher at Osaka Central
Advanced Mathematical Institute (OCAMI researcher), which is supported by MEXT Promotion of Distinctive Joint Research Center Program JPMXP0723833165. He is not employed there. We also thank them. The author would also like to thank the conference "Singularity theory of differentiable maps and its applications" (https://www.fit.ac.jp/$\sim$fukunaga/conf/sing202412.html) for an opportunity to present \cite{kitazawa2, kitazawa3, kitazawa4}. Comments presented there have motivated the author to study further including the present study. This conference is also supported by the Research Institute for Mathematical Sciences, an International Joint Usage/Research Center located in Kyoto University. \\
\ \\
{\bf Data availability.} \\
Data essentially related to our present study are all in the present file.


\begin{thebibliography}{25}
%	\bibitem{buchstaberpanov} V. M. Buchstaber and T. E. Panov, \textsl{Toric topology}, Mathematical Surveys and Monographs, Vol. 204, American Mathematical Society, Providence, RI, 2015.
%	\bibitem{burletderham} O. Burlet and G. de Rham, \textsl{Sur certaines applications g\'en\'eriques d'une vari\'et\'e close a $3$ dimensions dans le plan}, Enseign. Math. 20 (1974). 275--292.
	%	\bibitem{calabi} E. Calabi, Quasi-surjective mappings and a generation of Morse theory, Proc. U.S.-Japan Seminar in Differential Geometry, Kyoto, 1965, pp. 13--16.
	%
	%		\bibitem{cavicchioli} A. Cavicchioli, \textsl{Covering numbers of manifolds and critical points of a Morse function}, Israel. J. Math. 70 (1990), 279--304.
	% \bibitem{cerf} J. Cerf, \textsl{La stratification naturelle des espaces de fonctions deff\'erentiables r'eelles et le th'eor`eme de la pseudo-isotopie}, Inst. Hautes Etudes Sci. Publ. Math. 39 (1970), 5--173.
	%		\bibitem{choimasudasuh} S. Choi, M. Masuda and D. Y. Suh, \textsl{Topological classification of generalized Bott towers}, Trans. Amer. Math. Soc. 362 (2010), 1097--1112.
	%		\bibitem{cornealuptonopreatanre} O. Cornea, G. Lupton, J. Oprea and D. Tanr\'e, \textsl{Lusternik-Schnirelmann category}, Mathematical Surveys and Monographs, 103, Amer. Math. Soc., Providence, RI, 2003.
	%\bibitem{crowleyescher} D. Crowley and C. Escher, \textsl{A classification of $S^3$-bundles over $S^4$}, Differential. Geom. Appl. 18 (2003), 363--380, arXiv:0004147.
	%\bibitem{crowleynordstrom} D. Crowley and J. Nordstr\"{o}m, \textsl{The classification of $2$-connected $7$-manifolds}, Proc. London. Math. Soc. 119 (2019), 1--54, arXiv:1406.2226.

	\bibitem{bochnakcosteroy} J. Bochnak, M. Coste and M.-F. Roy, \textsl{Real algebraic geometry}, Ergebnisse der Mathematik und ihrer Grenzgebiete (3) [Results in Mathematics and Related Areas (3)], vol. 36, Springer-Verlag, Berlin, 1998. Translated from the 1987 French original; Revised by the authors.
%		\bibitem{bochnakkucharz} J. Bochnak and W. Kucharz, \textsl{Algebraic approximation of mappings into spheres}, Michigan Mathematical Journal, vol. 34, no. 1, 1987.
	\bibitem{bodinpopescupampusorea} A. Bodin, P. Popescu-Pampu and M. S. Sorea, \textsl{Poincar\'e-Reeb graphs of real algebraic domains}, Revista Matem\'atica Complutense, https://link.springer.com/article/10.1007/s13163-023-00469-y, 2023, arXiv:2207.06871v2.
%\bibitem{bott} R. Bott, \textsl{Nondegenerate critical manifolds}, Ann. of Math. 60 (1954), 248--261.
%\bibitem{carmesinschulz} S. Carmesin and A. Schulz, \textsl{Arrangements of orthogonal circles with many intersections}, Graph Drawing and Network Visualization (a conference paper), SPRINGER NATURE Link, 2021, arXiv:2106.03557v2. 
%\bibitem{costantino}  F. Costantino, \textsl{A short introduction to shadows of $4$-manifolds}, Fundamenta Mathematicae 251 no. 2 (2005), 427--442.
%\bibitem{costantinothurston} F. Costantino, D. Thurston, \textsl{$3$-manifolds efficiently bound $4$-manifolds}, J. Topol. 1 (2008),
%703--745.
%	\bibitem{delzant} T. Delzant, \textsl{Hamiltoniens p\'eriodiques et images convexes de l'application moment}, Bull. Soc. Math. France 116 (1988), No. 3, 315--339.
%\bibitem{ehresmann} C. Ehresmann, \textsl{Les connexions infinitesimales dans un espace fibre differentiable}, Colloque de Topologie, Bruxelles (1950), 29--55.
\bibitem{elredge} N. Elredge, \textsl{{\it Answer to} On finding polynomials that approximate a function and its derivative}, StackExchange, question 555712 (2013), https://math.stackexchange.com/questions/555712/on-finding-polynomials-that-approximate-a-function-and-its-derivative-extension.
%\bibitem{fujitakitabeppumitsuishi} H. Fujita, Y Kitabeppu and A. Mitsuishi, \textsl{Distance functions and convex bodies and symplectic toric manifolds}, arXiv:2003.02293.
%\bibitem{gelbukh} I. Gelbukh, \textsl{Loops in Reeb graphs of $n$-manifolds}, diskrete \& Computational Geometry, 59 (4) (2018), 843--863. 
%%\bibitem{gelbukh2} I. Gelbukh, \textsl{Approximation of Metric Spaces by Reeb Graphs: Cycle Rank of a Reeb Graph, the Co-rank of the Fundamental Group, and Large Components of Level Sets on Riemannian Manifolds}, Filomat (in press), arxiv:1903.00777.
\bibitem{gelbukh1} I. Gelbukh, \textsl{A finite graph is homeomorphic to the Reeb graph of a Morse-Bott function}, Mathematica Slovaca, 71 (3), 757--772, 2021; doi: 10.1515/ms-2021-0018. 
\bibitem{gelbukh2} I. Gelbukh, \textsl{Morse-Bott functions with two critical values on a surface}, Czechoslovak Mathematical Journal, 71 (3), 865--880, 2021; doi: 10.21136/CMJ.2021.0125-20. 
\bibitem{golubitskyguillemin} M. Golubitsky and V. Guillemin, \textsl{Stable Mappings and Their Singularities}, Graduate Texts in Mathematics (14), Springer-Verlag (1974).
%\bibitem{hempel} J. Hempel, \textsl{3- Manifolds}, AMS Chelsea Publishing, 2004. 
%\bibitem{hiratukasaeki} J. T. Hiratuka and O. Saeki, \textsl{Triangulating Stein factorizations of generic maps and Euler Characteristic formulas}, RIMS Kokyuroku Bessatsu B38 (2013), 61--89. 
%\bibitem{hiratukasaeki2} J. T. Hiratuka and O. Saeki, \textsl{Connected components of regular fibers of differentiable maps}, in "Topics on Real and Complex Singularities", Proceedings of the 4th Japanese-Australian Workshop (JARCS4), Kobe 2011,  World Scientific, 2014, 61--73. 
\bibitem{hirsch} M .W. Hirsch, \textsl{Smooth regular neighborhoods}, Ann. of Math., 76 (1962), 524--530.
%\bibitem{ishikawakoda} M. Ishikawa and Y. Koda, \textsl{Stable maps and branched shadows of $3$-manifolds}, Mathematische Annalen 367 (2017), no. 3, 1819--1863, arXiv:1403.0596.
%\bibitem{kitazawa1} N. Kitazawa, \textsl{On round fold maps} (in Japanese), RIMS Kokyuroku Bessatsu B38 (2013), 45--59.
%\bibitem{kitazawa2} N. Kitazawa, \textsl{On manifolds admitting fold maps with singular value sets of concentric spheres}, Doctoral Dissertation, Tokyo Institute of Technology (2014).
%\bibitem{kitazawa3} N. Kitazawa, \textsl{Fold maps with singular value sets of concentric spheres}, Hokkaido Mathematical Journal Vol.43, No.3 (2014), 327--359.
\bibitem{kitazawa1} N. Kitazawa, \textsl{On Reeb graphs induced from smooth functions on $3$-dimensional closed orientable manifolds with finitely many singular values}, Topol. Methods in Nonlinear Anal. Vol. 59 No. 2B, 897--912, arXiv:1902.08841.
%\bibitem{kitazawa1} N. Kitazawa, \textsl{On Reeb graphs induced from smooth functions on closed or open surfaces}, Methods of Functional Analysis and Topology Vol. 28 No. 2 (2022), 127--143, arXiv:1908.04340.
\bibitem{kitazawa2} N. Kitazawa, \textsl{Real algebraic functions on closed manifolds whose Reeb graphs are given graphs}, Methods of Functional Analysis and Topology Vol. 28 No. 4 (2022), 302--308, arXiv:2302.02339, 2023.
%\bibitem{kitazawa6} N. Kitazawa, \textsl{Explicit construction of explicit real algebraic functions and real algebraic manifolds via Reeb graphs}, Algebraic and geometric methods of analysis 2023 “The book of abstracts”, 49—51, this is the abstract book of the conference "Algebraic and geometric methods of analysis 2023" and published after a short review (https://www.imath.kiev.ua/$\sim$topology/conf/agma2023/), https://imath.kiev.ua/$\sim$topology/conf/agma2023/contents/abstracts/texts/kitazawa/kitazawa.pdf, 2023.
%\bibitem{kitazawa5} N. Kitazawa, \textsl{Notes on explicit special generic maps into Euclidean spaces whose dimensions are greater than $4$}, a revised version is submitted based on positive comments (major revision) by referees and editors after the first submission to a refereed journal, arXiv:2010.10078.

%\bibitem{kitazawa6} N. Kitazawa, \textsl{Round fold maps and the topologies and the differentiable structures of manifolds admitting explicit ones}, submitted to a refereed journal, arXiv:1304.0618.
%\bibitem{kitazawa0.5} N. Kitazawa, \textsl{Constructing fold maps by surgery operations and homological information of their Reeb spaces}, submitted to a refereed journal, arxiv:1508.05630.
%\bibitem{kitazawa0.6} N. Kitazawa, \textsl{Notes on fold maps obtained by surgery operations and algebraic information of their Reeb spaces}, arxiv:1811.04080.


%\bibitem{kitazawa6} N. Kitazawa, \textsl{On Reeb graphs induced from smooth functions on $3$-dimensional closed manifolds which may not be orientable}, a revised version is submitted to a refereed journal after based on positive comments by editors and referees after the second submission to a refreed journal, arXiv:2108.01300.
%\bibitem{kitazawa7} N. Kitazawa, \textsl{Realization problems of graphs as Reeb graphs of Morse functions with prescribed preimages}, submitted to a refereed journal, arXiv:2108.06913.
%\bibitem{kitazawa10} N. Kitazawa,\textsl{Round fold maps on $3$-dimensional manifolds and their integral and rational cohomology rings}, arXiv:2301.07008.
%\bibitem{kitazawa6} N. Kitazawa, \textsl{A class of naturally generalized special generic maps}, arXiv:2212.03174.
%\bibitem{kitazawa7} N. Kitazawa, \textsl{Construction of real algebraic functions with prescribed preimages}, submitted to a refereed journal as the second version based on positive comments by referees and editors, arXiv:2303.00953v3.
\bibitem{kitazawa3} N. Kitazawa, \textsl{Reconstructing real algebraic maps locally like moment-maps with prescribed images and compositions with the canonical projections to the $1$-dimensional real affine space}, the title has changed from previous versions, arXiv:2303.10723, 2024.


%\bibitem{kitazawa6} N. Kitazawa, \textsl{A note on real algebraic maps which are topologically special generic maps}, previous version(s) of the present article and the version arXiv:2303.00953v2 is submitted to a refereed journal, arXiv:2312.10646. 
\bibitem{kitazawa4} N. Kitazawa, \textsl{Some remarks on real algebraic maps which are topologically special generic maps}, submitted to a refereed journal, arXiv:2312.10646. 
\bibitem{kitazawa5} N. Kitazawa, \textsl{Arrangements of small circles for Morse-Bott functions and regions surrounded by them}, arXiv:2412.20626v3, 2025.
%\bibitem{kitazawa10} N. Kitazawa, \textsl{A note on cohomological structures of special generic maps}, a revised version is submitted based on positive comments by referees and editors after the third submission to a refereed journal.
%\bibitem{kitazawasaeki1} N. Kitazawa and O. Saeki, \textsl{Round fold maps on $3$-manifolds}, accepted for publication after a refereeing process and to appear in Algebraic \& Geometric Topology, arXiv:2105.00974.
		%	\bibitem{kitazawasaeki2} N. Kitazawa and O. Saeki, \textsl{Round fold maps of $n$-dimensional manifolds into ${\mathbb{R}}^{n-1}$}, submitted to a refereed journal, arXiv:2111.13510.
%\bibitem{ishikawakoda} M. Ishikawa, Y. Koda, \textsl{Stable maps and branched shadows of $3$-manifolds}, arXiv:1403.0596.
%\bibitem{kobayashisaeki} M. Kobayashi and O. Saeki, \textsl{Simplifying stable mappings into the plane from a global viewpoint}, Trans. Amer. Math. Soc. 348 (1996), 2607--2636.
\bibitem{kohnpieneranestadrydellshapirosinnsoreatelen} K. Kohn, R. Piene, K. Ranestad, F. Rydell, B. Shapiro, R. Sinn, M-S. Sorea and S. Telen, \textsl{Adjoints and Canonical Forms of Polypols}, to appear in Documenta Mathematica, arXiv:2108.11747.
\bibitem{kollar} J. Koll\'ar, \textsl{Nash's work in algebraic geometry}, Bulletin (New Series) of the American Mathematical Society (2) 54, 2017, 307--324.
%\bibitem{kucharz} W. Kucharz, \textsl{Some open questions in real algebraic geometry}, Proyecciones Journal of Mathematics, Vol. 41 No. 2 (2022), Universidad Cat\'olica del Norte Antofagasta, Chile, 437--448.
\bibitem{lellis} Camillo De Lellis, \textsl{The Masterpieces of John Forbes Nash Jr.}, The Abel Prize 2013--2017 (Helge Holden and Ragni Piene, eds.), Springer International Publishing, Cham, 2019, 391--499, https://www.math.ias.edu/delellis/sites/math.ias.edu.delellis/files/Nash\_Abel\_75.pdf, arXiv:1606.02551.
%\bibitem{martinezalfaromezasarmientooliveira} J. Martinez-Alfaro, I. S. Meza-Sarmiento and R. Oliveira, \textsl{Topological classification of simple Morse Bott functions on surfaces}, Contemp. Math. 675 (2016), 165--179.%
\bibitem{marzantowiczmichalak} W. Marzantowicz and L. P. Michalak, \textsl{Relations between Reeb graphs, systems of hypersurfaces and epimorphisms onto free groups}, Fund. Math., 265 (2), 97--140, 2024.
\bibitem{masumotosaeki} Y. Masumoto and O. Saeki, \textsl{A smooth function on a manifold with given Reeb graph}, Kyushu J. Math. 65 (2011), 75--84.
%\bibitem{maciasvirgospereirasaez} E. Mac\'ias-Virg\'os and M. J. Pereira-S\'aez, Height functions on compact symmetric spaces, Monatshefte f\"ur Mathematik 177 (2015), 119--140. 
\bibitem{michalak1} L. P. Michalak, \textsl{Realization of a graph as the Reeb graph of a Morse function on a manifold}. Topol. Methods in Nonlinear Anal. 52 (2) (2018), 749--762, arXiv:1805.06727.
\bibitem{michalak2} L. P. Michalak, \textsl{Combinatorial modifications of Reeb graphs and the realization problem}, Discrete Comput. Geom. 65 (2021), 1038--1060, arXiv:1811.08031.
%\bibitem{milnor} J. Milnor, \textsl{Singular points of complex hypersurfacs}, Annals of Mathematics Studies, No. 61, Princeton University Press, Princeton, N. J.; University of Tokyo Press, Tokyo, 1968.
%\bibitem{milnor} J. Milnor, \textsl{Lectures on the h-cobordism theorem}, Math. Notes, Princeton Univ. Press, Princeton, N.J. 1965.
%\bibitem{moise} E. E. Moise, \textsl{Affine Structures in $3$-Manifold{\rm :} V. The Triangulation Theorem and Hauptvermutung}, Ann. of Math., Second Series, Vol. 56, No. 1 (1952), 96--114.
%\bibitem{morin} B. Morin, \textsl{Formes canoniques des singulariti\'{e}s d\'{}une application diff\'{e}rentiable}, C. E. Acad. Sci. Paris 260 (1965), 5662--5665, 6503--6506.
%\bibitem{nash} J. Nash, \textsl{Real algbraic manifolds}, Ann. of Math. (2) 56 (1952), 405--421.
%\bibitem{ranicki} A. Ranicki, \textsl{Algebraic and geometric surgery}, https://www.maths.ed.ac.uk/~v1ranick/books/surgery.pdf, 2002.
%\bibitem{ramanujam} S. Ramanujam, \textsl{Morse theory of certain symmetric spaces}, J. Diff. Geom. 3 (1969), 213--229.
\bibitem{reeb} G. Reeb, \textsl{Sur les points singuliers d\'{}une forme de Pfaff compl\'{e}tement int\`{e}grable ou d\'{}une fonction num\'{e}rique}, Comptes Rendus
 Hebdomadaires des S\'{e}ances de I\'{}Acad\'{e}mie des Sciences 222 (1946), 847--849.
%\bibitem{saeki1} O. Saeki, \textsl{Notes on the topology of folds}, J. Math. Soc. Japan Volume 44, Number 3 (1992), 551--566.
%\bibitem{saeki1} O. Saeki, \textsl{Topology of special generic maps of manifolds into Euclidean spaces}, Topology Appl. 49 (1993), 265--293.
%\bibitem{saeki0.2} O. Saeki, \textsl{Topology of singular fibers of differentiable maps}, Lecture Notes in Math., Vol. 1854, Springer-Verlag, 2004. 
%\bibitem{saeki4} O. Saeki, \textsl{Morse functions with sphere fibers}, Hiroshima Math. J. Volume 36, Number 1 (2006),  141--170.
\bibitem{saeki1} O. Saeki, \textsl{Reeb spaces of smooth functions on manifolds}, International Mathematics Research Notices, maa301, Volume 2022, Issue 11, June 2022, 3740--3768, https://doi.org/10.1093/imrn/maa301.

\bibitem{saeki2} O. Saeki, \textsl{Reeb spaces of smooth functions on manifolds II}, Res. Math. Sci. 11, article number 24 (2024), https://link.springer.com/article/10.1007/s40687-024-00436-z.
%\bibitem{saekitakase} O. Saeki and M. Takase, \textsl{Desingularizing special generic maps}, Journal of G\"okova Geometry Topology (2013), 1--24.
%\bibitem{sakurai} S. Sakurai, Master Thesis, Kyushu. Univ..
% \bibitem{saekitakase} O. Saeki and M. Takase, \textsl{Desingularizing special generic maps}, Journal of Gokova Geometry Topology 7 (2013), 1--24.
%\bibitem{saeki2} O. Saeki, \textsl{Topology of special generic maps of manifolds into Euclidean spaces}, Topology Appl. 49 (1993), 265--293.
%\bibitem{saeki4} O. Saeki, \textsl{Singular fibers and $4$-dimensional cobordism group}, Pacific J. Math. 248 (2010), 233--256.
%\bibitem{saekisakuma} O. Saeki and K. Sakuma, \textsl{On special generic maps into ${\mathbb{R}}^3$}, Pacific J. Math. 184 (1998), 175--193.
%\bibitem{saekisuzuoka} O. Saeki and K. Suzuoka, \textsl{Generic smooth maps with sphere fibers} J. Math. Soc. Japan Volume 57, Number 3 (2005), 881--902.
\bibitem{sharko} V. Sharko, \textsl{About Kronrod-Reeb graph of a function on a manifold}, Methods of Functional Analysis and
 Topology 12 (2006), 389--396.
%\bibitem{shiota} M. Shiota, \textsl{Thom's conjecture on triangulations of maps}, Topology 39 (2000), 383--399.
\bibitem{sorea1} M. S. Sorea, \textsl{The shapes of level curves of real polynomials near strict local maxima},  Ph. D. Thesis, Universit\'e de Lille, Laboratoire Paul Painlev\'e, 2018.
\bibitem{sorea2} M. S. Sorea, \textsl{Measuring the local non-convexity of real algebraic curves}, Journal of Symbolic Computation 109 (2022), 482--509.
%\bibitem{sorea1} M. S. Sorea, \textsl{The shapes of level curves of real polynomials near strict local maxima},  Ph. D. Thesis, Universit\'e de %Lille, Laboratoire Paul Painlev\'e, 2018.
%\bibitem{sorea2} M. S. Sorea, \textsl{Measuring the local non-convexity of real algebraic curves}, J. Symbolic Compute. 109 (2022), 482--509.
%\bibitem{stong} R. E. Stong, \textsl{Notes on cobordsm theory}, Princeton Universty Press, 1968.
%\bibitem{takeuchi} M. Takeuchi, \textsl{Nice functions on symmetric spaces}, Osaka. J. Mat. (2) Vol. 6 (1969), 283--289.
%\bibitem{tamaki1} D. Tamaki, Algebraic Topology A Guide to literature,  http://pantodon.jp/index.rb?body=about, 2023. 
%\bibitem{tamaki2} D. Tamaki, Algebraic Topology A Guide to literature (Submanifold arrangement), http://pantodon.jp/index.rb?body=submanifold\_arrangement, 2023.
%\bibitem{tamaki3} D. Tamaki, Algebraic Topology A Guide to literature (Arrangement variations), http://pantodon.jp/index.rb?body=arrangement\_variations, 2023.

%\bibitem{thom} R. Thom, \textsl{Les singularites des applications differentiables}, Ann. Inst. Fourier (Grenoble) 6 (1955-56), 43--87.
%\bibitem{tognoli} A. Tognoli, \textsl{Su una congettura di Nash}, Ann. Scuola Norm. Sup. Pisa (3) 27 (1973), 167--185.
%\bibitem{turaev} Vladimir G. Turaev, \textsl{Topology of shadows}, Preprint, 1991.
%\bibitem{wall} C. T. C Wall, \textsl{Classification problems in differential topology -- {\rm I:} Classificationon handlebodies}, Topology 2 (1963), 253--261.
%\bibitem{wall2} C. T. C. Wall \textsl{Classification problems in differential topology -- {\rm II:} Diffeomorphismsof handlebodies}, Topology 2 (1963), 263--272.
%\bibitem{wall3} C. T. C. Wall, \textsl{Classification problems in differential topology -- {\rm Q:} Quadratic forms on finite groups and related topics}, Topology 2 (1963), 281--298.
%\bibitem{wall4} C. T. C. Wall, \textsl{Classification problems in differential topology -- {\rm III:} Applications to special cases}, Topology 3 (1965), 291--304.
%%\bibitem{wall5} C. T. C. Wall, \textsl{Classification problems in differential topology -- {\rm IV:} Thickenings}, Topology 5 (1966), 73--94.
%\bibitem{wall6} C. T. C. Wall, \textsl{Classification problems in differential topology -- {\rm VI:} Classification of |{\rm (}$s-1${\rm )}-connected {\rm (}$2s+1${\rm )}-manifolds}, Topology 6 (3) (1967), 273--296.
%\bibitem{whitney} H.  Whitney,  \textsl{On singularities of mappings of Euclidean spaces: I,  mappings of the plane into the plane},  Ann.  of Math.  62 (1955),  374--410. 

	%\bibitem{zhubr1} A. V. Zhubr, Closed simply-connected six-dimensional manifolds: proofs of classification theorems, Algebra i Analiz 12 (2000), no. 4, 126--230.
%\bibitem{zhubr2} A. V. Zhubr (responsible for the page), http://www.map.mpim-bonn.mpg.de/6-manifolds:\_1-connected.
\end{thebibliography}
\end{document}